\documentclass[12pt]{article}
\usepackage{amsmath,amssymb}
\usepackage[english]{babel}
\title{Multi-dimensional vector product}
\author{ Z.~K.~Silagadze
\vspace*{3mm} \\
Budker Institute of Nuclear Physics,  630 090,
Novosibirsk, Russia }
\date{}
\begin{document}
\maketitle

\begin{abstract}
It is shown that multi-dimensional generalization of the vector product is
only possible in seven dimensional space.
\end{abstract}


\maketitle

 
The three-dimensional vector product proved to be useful in various physical
problems. A natural question is whether multi-dimensional generalization 
of the vector product is possible. This apparently simple question has
somewhat unexpected answer, not widely known in physics community, that
generalization is only possible in seven dimensional space. In 
mathematics this fact was known since forties \cite{1}, but only recently
quite simple proof (in comparison to previous ones) was given by Markus
Rost \cite{2}. Below I present a version of this proof to make it
more accessible to physicists.

For contemporary physics seven-dimensional vector product represents not only
an academic interest. It turned out that the corresponding construction 
is useful in considering self-dual Yang-Mills fields depending only upon time 
(Nahm equations) which by themselves originate in the context of M-theory
\cite{3,4}. Other possible applications include Kaluza-Klein compactifications
of $d=11,\;N=1$ Supergravity \cite{5}. That is why I think that this beautiful
mathematical result should be known by a general audience of physicists.

Let us consider n-dimensional vector space $\mathbb{R}^n$ over the real
numbers with the standard Euclidean scalar product. Which properties we want 
the multi-dimensional bilinear vector product in 
$\mathbb{R}^n$ to satisfy? It is natural to choose as defining axioms the
following (intuitively most evident) properties of the usual three-dimensional
vector product:
\begin{equation}   
\vec{A}\times\vec{A}=0,
\label{eq1} \end{equation}
\begin{equation} 
(\vec{A}\times\vec{B})\cdot \vec{A}=(\vec{A}\times\vec{B})\cdot \vec{B}=
0, \label{eq2} \end{equation} 
\begin{equation}
|\vec{A}\times\vec{B}|=|\vec{A}|~|\vec{B}|,\;\; {\rm if} \; \; 
\vec{A}\cdot\vec{B}
=0. \label{eq3} \end{equation}
Here $|\vec{A}|^2=\vec{A}\cdot\vec{A}$ is the norm of the vector $\vec{A}$.

Then $$0=(\vec{A}+\vec{B})\times
(\vec{A}+\vec{B})=\vec{A}\times\vec{B}+\vec{B}\times\vec{A}$$ shows that
the vector product is anti-commutative. By the same trick one can prove
that $(\vec{A}\times\vec{B})\cdot\vec{C}$ is alternating in 
$\vec{A},\vec{B},\vec{C}$. For example $$0=((\vec{A}+\vec{C})\times\vec{B})
\cdot (\vec{A}+\vec{C})=(\vec{C}\times\vec{B})\cdot\vec{A}+
(\vec{A}\times\vec{B})\cdot\vec{C}$$ shows that 
$(\vec{C}\times\vec{B})\cdot\vec{A}=-(\vec{A}\times\vec{B})\cdot\vec{C}$.

For any two
vectors $\vec{A}$ and $\vec{B}$ the norm $|\vec{A}\times\vec{B}|^2$ 
equals to
$$\left | \left (\vec{A}-\frac{\vec{A}\cdot\vec{B}}
{|\vec{B}|^2}~\vec{B}\right)\times\vec{B}\right |^2=
\left | \vec{A}-\frac{\vec{A}\cdot\vec{B}}{|\vec{B}|^2}~\vec{B}
\right |^2|\vec{B}|^2=|\vec{A}|^2 |\vec{B}|^2-(\vec{A}\cdot\vec{B})^2.$$
Therefore for any two vectors we should have
\begin{equation}
(\vec{A}\times\vec{B})\cdot (\vec{A}\times\vec{B})=(\vec{A}\cdot\vec{A})
(\vec{B}\cdot\vec{B})-(\vec{A}\cdot\vec{B})^2.
\label{eq4} \end{equation}
Now consider
$$|\vec{A}\times (\vec{B}\times\vec{A})-(\vec{A}\cdot\vec{A})\vec{B}+
(\vec{A}\cdot\vec{B})\vec{A}|^2=$$ 
$$=|\vec{A}\times (\vec{B}\times\vec{A})|^2+
|\vec{A}|^4|\vec{B}|^2-(\vec{A}\cdot\vec{B})^2|\vec{A}|^2-
2|\vec{A}|^2(\vec{A}\times (\vec{B}\times\vec{A}))\cdot\vec{B}.$$
But this is zero because
$$|\vec{A}\times (\vec{B}\times\vec{A})|^2=|\vec{A}|^2|\vec{B}\times\vec{A}|^2
=|\vec{A}|^4|\vec{B}|^2-(\vec{A}\cdot\vec{B})^2|\vec{A}|^2$$
and
$$(\vec{A}\times (\vec{B}\times\vec{A}))\cdot\vec{B}=(\vec{B}\times\vec{A})
\cdot (\vec{B}\times\vec{A})=|\vec{A}|^2|\vec{B}|^2-(\vec{A}\cdot\vec{B})^2.$$ 
Therefore we have proved the identity
\begin{equation}
\vec{A}\times (\vec{B}\times\vec{A})=(\vec{A}\cdot\vec{A})\vec{B}-
(\vec{A}\cdot\vec{B})\vec{A}\;.
\label{eq5} \end{equation}
Note that arrangement of the brackets in the l.f.s. is in fact irrelevant
because the vector product is anti-commutative.

However, familiar identity
\begin{equation}
\vec{A}\times (\vec{B}\times\vec{C})=\vec{B}(\vec{A}\cdot\vec{C})-
\vec{C}(\vec{A}\cdot\vec{B})
\label{eq6} \end{equation}
does not follow in general from the intuitively evident properties
(\ref{eq1}-\ref{eq3}) of the vector product \cite{2}. To show this, let
us  introduce the ternary product \cite{6} (which is zero if (\ref{eq6}) is 
valid)
$$\{\vec{A},\vec{B},\vec{C}\}=\vec{A}\times (\vec{B}\times\vec{C})-
\vec{B}(\vec{A}\cdot\vec{C})+\vec{C}(\vec{A}\cdot\vec{B}).$$
Equation (\ref{eq5}) implies that this ternary product is alternating in its 
arguments. For example
$$0=\{\vec{A}+\vec{B},\vec{A}+\vec{B},\vec{C}\}=\{\vec{A},\vec{B},\vec{C}\}+
\{\vec{B},\vec{A},\vec{C}\}.$$
If $\vec{e}_i,\; i=1,\ldots, n$ is an orthonormal basis in the vector space, 
then
$$(\vec{e}_i\times\vec{A})\cdot(\vec{e}_i\times\vec{B})=
((\vec{e}_i\times\vec{B})\times\vec{e}_i)\cdot\vec{A}=
[\vec{B}-(\vec{B}\cdot\vec{e}_i)\vec{e}_i]\cdot\vec{A}$$
and, therefore,
\begin{equation}
\sum\limits_{i=1}^n(\vec{e}_i\times\vec{A})\cdot(\vec{e}_i\times\vec{B})=
(n-1)\vec{A}\cdot\vec{B}.
\label{eq7} \end{equation}
Using this identity we obtain
$$\sum\limits_{i=1}^n\{\vec{e}_i,\vec{A},\vec{B}\}\cdot 
\{\vec{e}_i,\vec{C},\vec{D}\}=$$ 
\begin{equation}
=(n-5)(\vec{A}\times \vec{B})\cdot
(\vec{C}\times \vec{D})+2(\vec{A}\cdot\vec{C})(\vec{B}\cdot\vec{D})-
2(\vec{A}\cdot\vec{D})(\vec{B}\cdot\vec{C}).
\label{eq8} \end{equation}
Hence
\begin{equation}
\sum\limits_{i,j=1}^n\{\vec{e}_i,\vec{e}_j,\vec{A}\}\cdot
\{\vec{e}_i,\vec{e}_j,\vec{B}\}=(n-1)(n-3)\vec{A}\cdot\vec{B}
\label{eq9} \end{equation}
and \cite{6}
\begin{equation}
\sum\limits_{i,j,k=1}^n\{\vec{e}_i,\vec{e}_j,\vec{e}_k\}\cdot
\{\vec{e}_i,\vec{e}_j,\vec{e}_k\}=n(n-1)(n-3).
\label{eq10} \end{equation}
The last equation shows that there exists 
some  $\{\vec{e}_i,\vec{e}_j,\vec{e}_k\}$ that 
is not zero if $n>3$. So equation (\ref{eq6}) is valid only for the usual 
three-dimensional vector product ($n=1$ case is, of course, not interesting 
because it corresponds to identically vanishing vector product). Surprisingly,
we do not have much choice for $n$ even in case when validity of 
(\ref{eq6}) is not required. In fact the space dimension $n$ should 
satisfy \cite{2} (see also \cite{Nieto})
\begin{equation}
n(n-1)(n-3)(n-7)=0.
\label{eq11}\end{equation} 
To prove this statement, let us note that using
$$\vec{A}\times (\vec{B}\times \vec{C})+(\vec{A}\times \vec{B})\times \vec{C}
=(\vec{A}+\vec{C})\times\vec{B}\times(\vec{A}+\vec{C})-
\vec{A}\times\vec{B}\times\vec{A}-\vec{C}\times\vec{B}\times\vec{C}=$$
$$=2\vec{A}\cdot\vec{C}~\vec{B}-\vec{A}\cdot\vec{B}~\vec{C}-
\vec{B}\cdot\vec{C}~\vec{A}$$
and
$$\vec{A}\times(\vec{B}\times(\vec{C}\times\vec{D}))=
\left . \frac{1}{2} \right [
\vec{A}\times (\vec{B}\times(\vec{C}\times\vec{D}))+
(\vec{A}\times \vec{B})\times(\vec{C}\times\vec{D})-$$
$$-(\vec{A}\times \vec{B})\times(\vec{C}\times\vec{D})-
((\vec{A}\times \vec{B})\times \vec{C})\times\vec{D}+
((\vec{A}\times \vec{B})\times\vec{C})\times\vec{D}+$$
$$+(\vec{A}\times (\vec{B}\times\vec{C}))\times\vec{D}-
(\vec{A}\times (\vec{B}\times\vec{C}))\times\vec{D}-
\vec{A}\times ((\vec{B}\times \vec{C})\times\vec{D})+ $$ 
$$+\vec{A}\times ((\vec{B}\times \vec{C})\times\vec{D})+
\vec{A}\times (\vec{B}\times ( \vec{C}\times\vec{D})) \left . \frac{}{}
\right ]$$
we can check the equation
$$\vec{A}\times \{\vec{B},\vec{C},\vec{D}\}=-\{\vec{A},\vec{B},\vec{C}
\times \vec{D}\}+\vec{A}\times(\vec{B}\times(\vec{C}\times\vec{D}))-
\{\vec{A},\vec{C},\vec{D}\times \vec{B} \}+ $$ $$+
\vec{A}\times(\vec{C}\times(\vec{D}\times\vec{B}))-
\{\vec{A},\vec{D},\vec{B}\times \vec{C} \}+
\vec{A}\times(\vec{D}\times(\vec{B}\times\vec{C}))=$$
$$=-\{\vec{A},\vec{B},\vec{C}\times \vec{D} \}
-\{\vec{A},\vec{C},\vec{D}\times \vec{B} \}
-\{\vec{A},\vec{D},\vec{B}\times \vec{C} \}+
3\vec{A}\times\{\vec{B},\vec{C},\vec{D}\}.$$
The last step follows from
$$3\{\vec{B},\vec{C},\vec{D}\}=\{\vec{B},\vec{C},\vec{D}\}+
\{\vec{C},\vec{D},\vec{B}\}+\{\vec{D},\vec{B},\vec{C}\} =$$ $$=
\vec{B}\times(\vec{C}\times\vec{D})+
\vec{C}\times(\vec{D}\times\vec{B})+
\vec{D}\times(\vec{B}\times\vec{C}).$$
Therefore the ternary product satisfies an interesting identity
\begin{equation}
2~\vec{A}\times \{\vec{B},\vec{C},\vec{D}\}=\{\vec{A},\vec{B},\vec{C}
\times \vec{D}\}+\{\vec{A},\vec{C},\vec{D}\times \vec{B} \}+
\{\vec{A},\vec{D},\vec{B}\times \vec{C} \} \;\;\;
\label{eq12}\end{equation}
Hence we should have
$$4\sum\limits_{i,j,k,l=1}^n|\vec{e}_i\times \{\vec{e}_j,\vec{e}_k,
\vec{e}_l\}|^2=$$ $$=
\sum\limits_{i,j,k,l=1}^n|\{\vec{e}_i,\vec{e}_j,
\vec{e}_k\times \vec{e}_l\}+\{\vec{e}_i,\vec{e}_k,\vec{e}_l\times
\vec{e}_j\}+\{\vec{e}_i,\vec{e}_l,\vec{e}_j\times\vec{e}_k\}|^2.$$
L.h.s. is easily calculated by means of (\ref{eq7}) and (\ref{eq10}):
$$4\sum\limits_{i,j,k,l=1}^n|\vec{e}_i\times \{\vec{e}_j,\vec{e}_k,
\vec{e}_l\}|^2=4n(n-1)^2(n-3).$$
To calculate the r.h.s. the following identity is useful
\begin{equation}
\sum\limits_{i,j=1}^n\{\vec{e}_i,\vec{e}_j,\vec{A}\}\cdot
\{\vec{e}_i,\vec{e}_j\times\vec{B},\vec{C}\}=-(n-3)(n-6)\vec{A}\cdot
(\vec{B}\times\vec{C})
\label{eq13}\end{equation}
which follows from (\ref{eq8}) and from the identity 
$$\sum\limits_{i=1}^n(\vec{e}_i\times\vec{A})\cdot((\vec{e}_i\times\vec{B})
\times\vec{C})=$$ $$=\sum\limits_{i=1}^n(\vec{e}_i\times\vec{A})\cdot
[2\vec{e}_i\cdot\vec{C}~\vec{B}-\vec{B}\cdot\vec{C}~\vec{e}_i-
\vec{e}_i\cdot\vec{B}~\vec{C}-\vec{e}_i\times(\vec{B}\times\vec{C})]=$$
$$=-(n-4)\vec{A}\cdot(\vec{B}\times\vec{C}).$$
Now, with (\ref{eq9}) and (\ref{eq13}) at hands, it becomes an easy task to 
calculate
$$\sum\limits_{i,j,k,l=1}^n|\{\vec{e}_i,\vec{e}_j,
\vec{e}_k\times \vec{e}_l\}+\{\vec{e}_i,\vec{e}_k,\vec{e}_l\times
\vec{e}_j\}+\{\vec{e}_i,\vec{e}_l,\vec{e}_j\times\vec{e}_k\}|^2=$$
$$=3n(n-1)^2(n-3)+6n(n-1)(n-3)(n-6)=3n(n-1)(n-3)(3n-13).$$
Therefore we should have
$$4n(n-1)^2(n-3)=3n(n-1)(n-3)(3n-13).$$
But 
$$3n(n-1)(n-3)(3n-13)-4n(n-1)^2(n-3)=5n(n-1)(n-3)(n-7)$$
and hence (\ref{eq11}) follows.

As we see, the space dimension must equal to the magic number seven if 
unique generalization of the ordinary three-dimensional vector product
is possible.

So far we only have shown that seven-dimensional  vector product can exist 
in principle. What about its detailed realization? To answer this question,
it is useful to realize that the vector products are closely related to 
composition algebras \cite{1} (in fact these two notions are equivalent 
\cite{2}). Namely, for any composition algebra with unit element $e$
we can define the vector product in the subspace orthogonal to $e$ by
$x\times y=\frac{1}{2}(xy-yx)$. Therefore from a viewpoint of 
composition algebra, the
vector product is just the commutator divided by two. 
According to Hurwitz theorem \cite{7} the only composition
algebras are real numbers, complex numbers, quaternions and octonions. The
first two of them give identically zero vector products. Quaternions produce
the usual three-dimensional vector product. The seven-dimensional vector 
product is generated by octonions  \cite{8}. It is interesting to note that
this seven-dimensional vector product is covariant with respect to  
smallest exceptional Lie group $G_2$ \cite{9} which is the automorphism
group of octonions.

Using the octonion multiplication table \cite{8} one can realize the
seven-dimensional vector product as follows
\begin{equation}
\vec{e}_i\times\vec{e}_j=
\sum\limits_{k=1}^7 f_{ijk}\vec{e}_k,\;\;\;
i,j=1,2,\ldots,7,
\label{eq14}\end{equation}
where $f_{ijk}$ is totally antisymmetric $G_2$-invariant tensor and the
only nonzero independent components are
$$f_{123}=f_{246}=f_{435}=f_{651}=f_{572}=f_{714}=f_{367}=1.$$
Note that in contrast to the three-dimensional case
$f_{ijk}f_{kmn}\ne\delta_{im}\delta_{jn}-\delta_{in}\delta_{jm}.$ Instead
we have
\begin{equation}
f_{ijk}f_{kmn}=g_{ijmn}+\delta_{im}\delta_{jn}-\delta_{in}\delta_{jm}
\label{eq15} \end{equation}
where
$$g_{ijmn}=\vec{e}_i\cdot\{\vec{e}_j,\vec{e}_m,\vec{e}_n\}.$$
In fact $g_{ijmn}$ is totally antisymmetric $G_2$-invariant tensor.
For example
$$g_{ijmn}=\vec{e}_i\cdot\{\vec{e}_j,\vec{e}_m,\vec{e}_n\}=
-\vec{e}_i\cdot\{\vec{e}_m,\vec{e}_j,\vec{e}_n\}=$$ $$=
-\vec{e}_i\cdot(\vec{e}_m\times(\vec{e}_j\times\vec{e}_n))+(\vec{e}_i\cdot
\vec{e}_j)(\vec{e}_m\cdot\vec{e}_n)-(\vec{e}_i\cdot\vec{e}_n)
(\vec{e}_m\cdot\vec{e}_j)=$$ $$=-\vec{e}_m\cdot(\vec{e}_i\times
(\vec{e}_n\times\vec{e}_j))+(\vec{e}_i\cdot\vec{e}_j)(\vec{e}_m
\cdot\vec{e}_n)-(\vec{e}_i\cdot\vec{e}_n)(\vec{e}_m\cdot\vec{e}_j)=$$ $$=
-\vec{e}_m\cdot\{\vec{e}_i,\vec{e}_n,\vec{e}_j\}=
-\vec{e}_m\cdot\{\vec{e}_j,\vec{e}_i,\vec{e}_n\}=-g_{mjin}.$$
The only nonzero independent components are
$$g_{1254}=g_{1267}=g_{1364}=g_{1375}=g_{2347}=g_{2365}=g_{4576}=1.$$
In conclusion, generalization of the vector product we have considered
is only possible in seven-dimensional space and
is closely related to octonions -- the largest composition algebra
which ties up together many exceptional structures in mathematics \cite{9}.
In a general case of p-fold vector products other options
arise \cite{1,5,10}. We recommend that interested reader consults  
references to explore these possibilities and possible physical applications.


\end{document}